\documentclass[11pt]{article}
\usepackage[margin=1in]{geometry}
\usepackage{amsmath,amssymb,amsthm}
\usepackage[hidelinks]{hyperref}
\usepackage[expansion=false]{microtype}

\newtheorem{theorem}{Theorem}
\newtheorem{proposition}[theorem]{Proposition}
\newtheorem{lemma}[theorem]{Lemma}
\newtheorem{corollary}[theorem]{Corollary}
\theoremstyle{remark}

\newcommand{\F}{\mathbb{F}}

\renewcommand{\DH}{\mathrm{DH}}

\title{A Little Advice Suffices:\\ Diffie--Hellman Is as Hard as Discrete Logarithm}
\newif\ifanon
\anonfalse
\ifanon
\author{Anonymous}
\else
\author{Yuwen Zhou\\[2pt] \normalsize University of British Columbia Okanagan\\
\normalsize\texttt{zhouyw07@student.ubc.ca}}
\fi
\date{}

\begin{document}
\maketitle

\begin{abstract}
Is computing a Diffie--Hellman key as hard as computing a discrete logarithm? Maurer showed that it is,
given a few bits of advice, provided that short intervals contain enough integers with only small prime
factors; that proviso has stood for thirty years. We remove it. For every prime $p$ there is an advice
string of $O(\log p)$ bits, depending only on $p$, with which discrete logarithms in every group of order $p$
are computed in polynomial time from a Diffie--Hellman oracle. The reduction is the classical one; what was
missing was a proof that short intervals contain polylogarithmically smooth integers. We show that a
zero-free half-plane for the Riemann zeta function, which is now known, gives the expected number of
$(\log x)^K$-smooth integers in every interval $[x,x+x^\theta]$ with $\theta>17/30$. With smoothness this
small the saddle point of the counting integral lies just to the left of the line $\operatorname{Re}s=1$,
and a zero-free strip of any fixed width absorbs the shift once $K$ is large. Genus-2 Jacobians, whose
orders fill such intervals, then turn one of these integers into the order of an auxiliary group. The same
input removes the smoothness assumption from Maurer and Raub's impossibility of field-homomorphic one-way
permutations and from the yes/no-query reduction of Koblitz, Menezes and Shparlinski. It also makes Wolf's
characterisation unconditional: discrete logarithm and Diffie--Hellman are generically equivalent exactly
when the group order has no large repeated prime factor.
\end{abstract}

\section{Introduction}\label{sec:intro}

The Diffie--Hellman protocol~\cite{DH76} is secure only if the \emph{computational Diffie--Hellman problem}
is hard: given $g^a$ and $g^b$ in a cyclic group, compute $g^{ab}$. What cryptographers actually believe
hard in the groups they deploy is the \emph{discrete logarithm problem}: given $g^a$, compute $a$. A
discrete-logarithm algorithm breaks Diffie--Hellman at once. Whether the converse holds, so that the
protocol is exactly as secure as discrete logarithms, has been asked since the protocol appeared.

The known answers all follow one idea, due to den Boer~\cite{dB88} and Maurer~\cite{Mau94}. A
Diffie--Hellman oracle lets one multiply hidden exponents, so the hidden $a$ can be treated as an element of
a field one can compute with but cannot read. Put $a$ into an auxiliary algebraic group over $\F_p$ whose
order has only small prime factors, solve the discrete logarithm there by Pohlig--Hellman, and read $a$ back.
Everything works once such an auxiliary group exists. Maurer used elliptic curves, whose orders fill the
interval $[p+1-2\sqrt p,\,p+1+2\sqrt p]$, and needed this interval to contain a polylogarithmically smooth
integer, which is still unknown. Maurer and Wolf~\cite{MW99,Wol99}, Boneh and Lipton~\cite{BL96}, and Maurer
and Raub~\cite{MR07} built on the same assumption.

Without it, the best reduction ran in time
$\exp((\log p)^{2/3+o(1)})$~\cite{Wol99,MW00}. Under the Riemann Hypothesis, Koblitz, Menezes and
Shparlinski~\cite{KMS10} reached $\exp(O(\sqrt{\log p\log\log p}))$ with the help of $\varepsilon\log p$
answers from an all-knowing oracle; they called a polynomial-time reduction a ``holy grail'' that ``seems a
long way off'', and named the missing input: smooth numbers in short intervals for smoothness bounds below
$\exp(\sqrt{\log x})$.

\paragraph{The missing input.}
Genus-2 Jacobians widen the window. Their orders over $\F_p$ fill
$[p^2-\tfrac12p^{3/2},\,p^2+\tfrac12p^{3/2}]$, an interval of the form $[x,x+x^{3/4}]$, and Lenstra, Pila and
Pomerance~\cite{LPP02} showed that all but at most $p$ integers there are such orders. Elliptic curves give
only intervals of the form $[x,x+4x^{1/2}]$. What is needed is
more than $p$ integers in that window whose prime factors are all at most $(\log p)^{K}$. Unconditional
methods give smooth numbers in $[x,x+x^\theta]$ only for smoothness bounds
$y\ge\exp((\log x)^{2/3+o(1)})$, because they rely on the Vinogradov--Korobov zero-free region. Younis
recently reached $y=(\log x)^K$ under the Riemann Hypothesis~\cite{You24}.

\paragraph{A zero-free half-plane suffices.}
We prove that any half-plane $\operatorname{Re}s>1-\delta$ free of zeros of $\zeta$ is enough
(Theorem~\ref{thm:smooth}): for every $\theta>17/30$ and every fixed $K$ above a threshold depending on
$\theta$ and $\delta$, the number of $y$-smooth integers in $[x,x+h]$, $y=(\log x)^K$, is asymptotic to the
count predicted by the density of $y$-smooth integers near $x$ whenever $x^\theta\le h=o(x)$. This range
contains the genus-2 window, $\theta=3/4$, but not the elliptic-curve window, $\theta=1/2$. Such a
half-plane, with $\delta=1/8$, is now known~\cite{OAI26}, so our results are unconditional. The count is a
contour integral of $\zeta(s,y)x^s$, which one shifts to the left of the saddle point $\alpha\approx1-1/K$. A
zero of $\zeta$ at distance $\nu$ from the line $\operatorname{Re}s=1$ costs a factor $x^{(1-\alpha)-\nu}$,
the half-plane guarantees $\nu\ge\delta$, and the zero-density estimate of Guth and Maynard~\cite{GM24}
bounds how many zeros there are. Polylogarithmic smoothness moves the saddle point only about $1/K$ from the
line, and the strip $1-\delta<\operatorname{Re}s<1$ absorbs this once $K$ is large in terms of $\delta$.
Earlier arguments placed the saddle point deep inside the zero-free region, which forced much larger
$y$~\cite{Zha26}; what lets it leave that region is a local shift bound that holds for every $y\ge2$
(Lemma~\ref{lem:shift} in Appendix~\ref{app:smooth}).

\paragraph{Results.}
With this input, Maurer's method goes through with genus-2 Jacobians.

\begin{theorem}\label{thm:main}
There are a constant $c$ and a probabilistic algorithm $\mathcal A$ with the following property. For every
prime $p$ there is an advice string $\sigma_p$ of at most $30\lceil\log_2p\rceil+12$ bits such that for every
cyclic group $G$ of order $p$, every generator $g$ and every $a$, the algorithm $\mathcal A$ given $p$,
$\sigma_p$, $g^a$ and a Diffie--Hellman oracle for $(G,g)$ outputs $a$ in expected time $O((\log p)^c)$.
\end{theorem}

The advice depends only on $p$: one string serves every group of order $p$, for instance every elliptic
curve group of that order, and every generator. Sections~\ref{sec:fewer} and~\ref{sec:consequences} derive
the following.
\begin{itemize}
\item \emph{Imperfect adversaries.} An algorithm that solves Diffie--Hellman on a fraction $\varepsilon$ of
  inputs yields discrete logarithms in time $\mathrm{poly}(\log p,1/\varepsilon)$
  (Corollary~\ref{cor:imperfect}).
\item \emph{A complete characterisation.} In cyclic groups of order $n$, the generic non-uniform complexity
  of computing discrete logarithms from a Diffie--Hellman oracle lies between $\sqrt q/2-1$ and
  $(\log n)^{O(1)}(1+\sqrt q)$, where $q$ is the largest prime with $q^2\mid n$
  (Theorem~\ref{thm:char}). Wolf proved this dichotomy under his smoothness assumption~\cite{Wol99}.
\item \emph{Black-box fields.} Extraction from a black-box field is polynomial-time with advice for every
  finite field. So, for characteristic fixed by the security parameter, there are no field-homomorphic
  one-way permutations against non-uniform adversaries (Corollary~\ref{cor:bbf}); Maurer and Raub proved
  this for small characteristic~\cite{MR07}.
\item \emph{Less advice.} With $\varepsilon\log p$ yes/no queries instead of advice the reduction is still
  polynomial (Theorem~\ref{thm:eps}); with a one-time precomputation of $p^{\varepsilon}$ steps per prime it
  costs $(\log p)^{O(1/\varepsilon)}$ per instance (Corollary~\ref{cor:precomp}); with no help at all it takes
  $\exp(O(\sqrt{\log p\log\log p}))$ (Theorem~\ref{thm:uniform}).
\end{itemize}

\paragraph{Contributions.}
\begin{itemize}
\item \emph{A theorem on smooth numbers.} A zero-free half-plane for $\zeta$ gives the expected count of
  polylogarithmically smooth integers in short intervals.
\item \emph{An unconditional equivalence.} Discrete logarithm and Diffie--Hellman are polynomially
  equivalent in every group of prime order, with logarithmic advice that depends only on the order, and the
  classical consequences of this equivalence hold without assumptions.
\item \emph{An explanation.} The reduction was never the obstacle; the existence of smooth group orders was.
  A zero-free strip of any fixed width provides them, because polylogarithmic smoothness keeps the saddle
  point close to the line $\operatorname{Re}s=1$.
\end{itemize}

\paragraph{Technical overview.}
The reduction (Section~\ref{sec:reduction}) is Maurer's, written for genus~2 as in~\cite{Wol99,KMS10}. A
Diffie--Hellman oracle turns the group into a black-box field: $g^a$ represents $a$, the group operation adds,
the oracle multiplies, and inversion is exponentiation by $p-2$. The advice names a genus-2 curve
$Y^2=f(X)$ whose Jacobian $J$ has $(\log p)^K$-smooth order, together with a basis of $J(\F_p)$. The
algorithm shifts $a$ by a random $r$ until $f(a+r)$ is a square, takes an implicit square root, and so holds
the point $(a+r,b)$ of the curve, hence the class $D=[(a+r,b)-\infty]\in J(\F_p)$, in encoded form. Cantor's
algorithm runs on encoded coordinates. A Pohlig--Hellman routine for a non-cyclic group with an explicit
basis and an encoded target writes $D$ in the basis; the explicit representation of $D$ then shows $a+r$.
Section~\ref{sec:smooth} and Appendix~\ref{app:smooth} prove the smooth-number theorem, from which
Section~\ref{sec:smooth} derives Proposition~\ref{prop:goodcurve}; this is the only place where analytic
number theory enters.

\paragraph{Related work.}
Den Boer~\cite{dB88} proved the equivalence for primes $p$ with $p-1$ smooth. Maurer~\cite{Mau94} and Maurer and
Wolf~\cite{MW96,MW99,MW00} introduced auxiliary groups, advice, and the amplification of faulty oracles, and
proved generic lower bounds for orders with repeated prime factors~\cite{MW98}; Wolf's thesis~\cite{Wol99}
collects these results and the genus-2 variant. Boneh and Lipton~\cite{BL96} introduced black-box fields and a
uniform subexponential algorithm under a smoothness assumption; Maurer and Raub~\cite{MR07} treated extension
fields. Koblitz, Menezes and Shparlinski~\cite{KMS10} studied reductions with auxiliary oracles, and May and
Schneider~\cite{MS23} gave concrete reductions for standardised curves. Lenstra, Pila and
Pomerance~\cite{LPP93,LPP02} introduced genus-2 Jacobians to obtain rigorous smoothness tests. On smooth numbers
in short intervals, Soundararajan~\cite{Sou10} and Matom\"aki and Radziwi\l\l~\cite{MR16} treat large smoothness
bounds, Younis~\cite{You24} the polylogarithmic range under the Riemann Hypothesis (our proof makes his
constant explicit, Appendix~\ref{app:smooth}), and Zhang~\cite{Zha26} a fixed zero-free strip for
$y=\exp(\sqrt{\log x\log\log x})$. Our reduction is the classical one; what is new is the proof that the
auxiliary groups exist.

\section{Preliminaries}\label{sec:prelim}

\paragraph{Groups and oracles.}
Throughout, $p$ is a prime and $\ell=\lceil\log_2 p\rceil$. A group $G$ of order $p$ with generator $g$ is
given by encodings of its elements together with algorithms (or oracles) for the group operation,
inversion and equality. The Diffie--Hellman oracle for $(G,g)$ maps $(g^a,g^b)$ to $g^{ab}$. The discrete
logarithm problem asks for $a$ given $g^a$. Running times count group operations, oracle calls and bit
operations. An \emph{advice string} for $p$ is a bit string that may depend on $p$ but not on the group, the
generator or the instance. We write $L_p(\tfrac12,c)=\exp\bigl((c+o(1))\sqrt{\log p\,\log\log p}\bigr)$.

\paragraph{Black-box fields.}
Following Boneh and Lipton~\cite{BL96} and Maurer and Raub~\cite{MR07}, a black-box field isomorphic to
$\F_p$ is an unknown injective encoding $s\colon\F_p\to\{0,1\}^m$ together with oracles that return
$s(a\pm b)$, $s(ab)$ and $s(a^{-1})$ (for $a\neq 0$) from $s(a)$ and $s(b)$; the encodings $s(0)$ and $s(1)$
are given, and equality of field elements is equality of encodings. The \emph{extraction problem} asks for
$a$ given $s(a)$. For an explicit $c\in\F_p$ one computes $s(c)$ from $s(1)$ with $O(\ell)$ additions.

\begin{lemma}[Diffie--Hellman oracles implement black-box fields~\cite{Mau94,BL96,MR07}]\label{lem:dh-bbf}
Let $G=\langle g\rangle$ have order $p$ and put $s(a)=g^a$. Then addition, subtraction and equality are the
group operation, inversion and equality of $G$; multiplication is one call to the Diffie--Hellman oracle;
and inversion, $s(a^{-1})=s(a^{p-2})$, takes at most $2\ell$ oracle calls. In particular, an algorithm that
solves extraction for black-box fields isomorphic to $\F_p$ with $T$ field operations solves the discrete
logarithm problem in every group of order $p$ with $O(\ell T)$ group operations and oracle calls.
\end{lemma}

An \emph{$\F_p$-program with equality tests} computes on field elements (explicit constants and one hidden
input $a$) with $+$, $-$, $\times$, and division by elements that it has itself tested to be non-zero, and it
may branch on equality tests between its data. Such a program runs unchanged when $a$ is given only as
$s(a)$: replace every datum $d$ by $s(d)$, every operation by an oracle call and every test by a comparison
of encodings. It then outputs the encodings of its results and all of its branch bits.
Horner evaluation, square-and-multiply, the Tonelli--Shanks square-root algorithm with an explicit
quadratic non-residue~\cite{Sha73}, and Cantor's algorithm for Jacobians of hyperelliptic
curves~\cite{Can87} are such programs.

\paragraph{Smooth numbers.}
An integer $n\ge1$ is \emph{$y$-smooth} if all its prime factors are at most $y$, and
$\Psi(x,y)=\#\{n\le x: n \text{ is } y\text{-smooth}\}$. We write $u=\log x/\log y$, and $\alpha=\alpha(x,y)$
for the saddle point, the $\sigma>0$ minimising $x^\sigma\zeta(\sigma,y)$ with
$\zeta(s,y)=\prod_{q\le y}(1-q^{-s})^{-1}$~\cite{HT86}. For fixed $K>1$ one has
$\Psi(x,(\log x)^K)=x^{1-1/K+o(1)}$~\cite{dB66,CEP83}.

\paragraph{Genus-2 Jacobians.}
For a squarefree quintic $f\in\F_p[X]$, the equation $Y^2=f(X)$ defines a curve $C_f$ of genus~2 with one
point $\infty$ at infinity; $J_f$ denotes its Jacobian. Every class in $J_f(\F_p)$ has a unique
\emph{reduced Mumford representation} $(u,v)$: $u,v\in\F_p[X]$, $u$ monic, $\deg v<\deg u\le2$, and
$u\mid v^2-f$~\cite{Can87}. An affine point $(a,b)\in C_f(\F_p)$ gives the class $[(a,b)-\infty]$ with
representation $(X-a,\,b)$. Cantor's algorithm adds two reduced representations with $O(1)$ field
operations, among them $O(1)$ inversions. The group $J_f(\F_p)$ has order in
$[(\sqrt p-1)^4,(\sqrt p+1)^4]$ and is a product of at most four cyclic groups, since for every prime $q$ the
$q$-torsion of $J_f(\overline{\F}_p)$ has rank at most~4~\cite{Mum70}. If $f$ has leading coefficient $c$,
then $\tilde f(X)=c^{-6}f(cX)$ is monic and $(X,Y)\mapsto(cX,c^3Y)$ is an isomorphism
$C_{\tilde f}\to C_f$ over $\F_p$, so we may take $f$ monic.

\section{Smooth numbers in short intervals}\label{sec:smooth}

\begin{theorem}\label{thm:smooth}
Let $0<\delta\le\tfrac12$ and suppose $\zeta(s)\neq0$ for $\operatorname{Re}s>1-\delta$. Let
$17/30<\theta\le1$ and
\[
K>K^*(\theta,\delta):=\frac{13+30\delta}{\delta\,(30\theta-17)} .
\]
There are $\varepsilon_0>0$, $C$ and $x_0$, depending only on $\delta,\theta,K$, such that for all $x\ge x_0$,
$x^\theta\le h\le x$ and $(\log x)^K\le y\le x^{\varepsilon_0}$,
\begin{equation}\label{eq:asym}
\left|\frac{\Psi(x+h,y)-\Psi(x,y)}{\alpha\,(h/x)\,\Psi(x,y)}-1\right|
\le C\left(\frac1u+\frac hx\cdot\frac{\log(u+1)}{\log y}\right),
\end{equation}
and in particular $\Psi(x+h,y)-\Psi(x,y)\ge\tfrac13\,h\,\Psi(x,y)/x$.
\end{theorem}

Theorem~\ref{thm:smooth} also gives the hyperelliptic curve method of Lenstra, Pila and Pomerance a proved
running time $L_p(\tfrac12,2)(\log n)^{O(1)}$ for splitting an integer $n$ whose least prime factor is
$p$~\cite{Comp26b}.

We apply Theorem~\ref{thm:smooth} with the half-plane of~\cite{OAI26}:

\begin{theorem}[\cite{OAI26}]\label{thm:qrh}
$\zeta(s)\neq0$ for $\operatorname{Re}s>7/8$.
\end{theorem}

Any fixed $\delta>0$ suffices for every result in this paper; only the smoothness exponent $K$, and hence the
degree of the polynomial in Theorem~\ref{thm:main}, depends on it.

\paragraph{How the proof goes.}
Perron's formula writes $\Psi(x+h,y)-\Psi(x,y)$ as $\frac1{2\pi i}\int\zeta(s,y)\frac{(x+h)^s-x^s}{s}\,ds$
over a vertical segment through the saddle point $\alpha$, with height $H=x^{1+1/K+\varepsilon_1}/h$. As in
Younis~\cite{You24} and Zhang~\cite{Zha26}, the segment is moved onto a contour that stays on
$\operatorname{Re}s=\alpha$ near the real axis, bends left along an exponential path, wraps around each zero
of $\zeta$ at a distance proportional to its distance from the line $\operatorname{Re}s=1$, and otherwise runs
along $\operatorname{Re}s=(1+\varepsilon)/2$ for a small fixed $\varepsilon>0$. The integrand is holomorphic
in $\operatorname{Re}s>0$. The central part gives the main term $\alpha(h/x)\Psi(x,y)$ by the saddle-point
method of Hildebrand and Tenenbaum~\cite{HT86}. Elsewhere one needs a saving over $x^\alpha\zeta(\alpha,y)$,
and it comes from moving left: as long as no zero is near, $|x^s\zeta(s,y)|$ drops by
$x^{-(1-\varepsilon)(\alpha-\sigma)}$ when $\operatorname{Re}s$ drops from $\alpha$ to $\sigma$
(Lemma~\ref{lem:shift}). This local shift holds for \emph{every} $y\ge2$. Earlier versions required
$y\ge\exp((\log\log x)^2)$ or a saddle point within $\varepsilon\delta$ of the line $\operatorname{Re}s=1$;
without these conditions $y$ can be polylogarithmic, and then the saddle point sits about $1/K$ from the
line. A zero at distance $\nu\ge\delta$ from the line costs $x^{(1-\alpha)-(1-\varepsilon)^2\nu}$, there are
at most $T^{(30/13)\nu}=x^{(30/13)(1-\theta+1/K)\nu}$ of them, and the sum over zeros is small once
$K>K^*(\theta,\delta)$. The full proof is in Appendix~\ref{app:smooth}.

\paragraph{From smooth numbers to good curves.}
\begin{theorem}[{Lenstra, Pila, Pomerance~\cite[Thm.~1.1]{LPP02}}]\label{thm:lpp}
Let $p>8100$ be prime. For all but at most $p$ integers $z$ in $[p^2-\tfrac12p^{3/2},\,p^2+\tfrac12p^{3/2}]$
there are at least $p^{9/2}/(24000\log(p^2)\log\log(p^2))$ quintic polynomials $f\in\F_p[X]$ with non-zero
discriminant such that $\#J_f(\F_p)=z$.
\end{theorem}

\begin{proposition}\label{prop:goodcurve}
Let $K_0=25$, fix $K\ge K_0$, and put $x=p^2-\tfrac12p^{3/2}$ and $y=(\log x)^K$. For all sufficiently large
primes $p$, the interval $(x,x+p^{3/2}]$ contains $S_p\ge p^{3/2-2/K-o(1)}>2p$ integers that are $y$-smooth,
and there is a monic quintic $f\in\F_p[X]$ with non-zero discriminant such that $\#J_f(\F_p)$ is $y$-smooth.
\end{proposition}

\begin{proof}
Apply Theorems~\ref{thm:smooth} and~\ref{thm:qrh} with $\delta=1/8$, $\theta=3/4$ and $h=p^{3/2}$; then
$x^{3/4}<h\le x$ because $x<p^2$, and $K>K^*=268/11$. By the estimate
$\Psi(x,(\log x)^K)=x^{1-1/K+o(1)}$ of Section~\ref{sec:prelim} and $x\le p^2$, Theorem~\ref{thm:smooth} gives
$S_p\ge\tfrac13p^{3/2}x^{-1/K-o(1)}\ge p^{3/2-2/K-o(1)}$, which exceeds $2p$ for large $p$ because
$K\ge K_0>4$. At most $p$ of these integers are exceptional in Theorem~\ref{thm:lpp}, so one of them is
$\#J_f(\F_p)$ for some quintic $f$, which can be made monic as in Section~\ref{sec:prelim}.
\end{proof}

We call a quintic $f\in\F_p[X]$, and its coefficient vector, \emph{good} if its discriminant is non-zero and
$\#J_f(\F_p)$ is $y$-smooth.

\section{The reduction}\label{sec:reduction}

This section proves Theorem~\ref{thm:main} from Proposition~\ref{prop:goodcurve}. Fix $K$ and
$y=(\log x)^K$ with $x=p^2-\tfrac12p^{3/2}$ as in Proposition~\ref{prop:goodcurve}; note
$y\le(2\ell)^K$. For the finitely many primes below the threshold of Proposition~\ref{prop:goodcurve} the
algorithm uses exhaustive search.

\subsection{The advice}
Let $f$ be a monic quintic with non-zero discriminant such that $J=J_f(\F_p)$ has $y$-smooth order $N$. Write
$J=\bigoplus_{i=1}^4\langle B_i\rangle$ with $B_i$ of order $n_i$ and $n_1\mid n_2\mid n_3\mid n_4$ (some
$n_i$ may be~1), and let $\xi$ be a quadratic non-residue modulo~$p$. The advice is
\[
\sigma_p=\bigl(f;\ (u_i,v_i)_{i\le4};\ n_1,\dots,n_4;\ \xi\bigr),
\]
where $(u_i,v_i)$ is the Mumford representation of $B_i$. It has at most
$5\ell+4(4\ell+2)+4(2\ell+1)+\ell=30\ell+12$ bits and depends only on~$p$. From it the algorithm computes
$N=n_1n_2n_3n_4$ and factors $N$ by trial division up to $y$.

\subsection{Step 1: putting the hidden element on the curve}
This is Maurer's embedding~\cite{Mau94}, as adapted to genus~2 by Wolf and by Koblitz, Menezes and
Shparlinski~\cite{Wol99,KMS10}. Given $s(a)$, pick $r\in\F_p$ at random, form $s(a')$ with $a'=a+r$, and
compute $s(w)$ for $w=f(a')$ and $s(w^{(p-1)/2})$. If $w=0$ put $b=0$; if $w^{(p-1)/2}=1$ compute $s(b)$ with
$b^2=w$ by Tonelli--Shanks with the explicit $\xi$; otherwise pick a new $r$. Now $D=(X-a',\,b)$, held
implicitly, is the reduced representation of $[(a',b)-\infty]\in J$.

The affine curve has $p-t$ points with $|t|\le4\sqrt p$ by Weil's bound, so the number of $r$ with
$f(a+r)$ a square or zero is at least $(p-4\sqrt p)/2\ge p/3$ for $p\ge150$. Hence at most three trials are
expected, each with $O(\ell)$ field operations, followed by one square root with $O(\ell^2)$.

\subsection{Step 2: Pohlig--Hellman with a hidden target}
The core of the reduction is the following routine. It is the textbook digit-by-digit method written for a
non-cyclic group whose basis is known explicitly while the target is known only implicitly.

\begin{lemma}\label{lem:ph}
Let $H=\bigoplus_{i=1}^m\langle C_i\rangle$ be a finite abelian $q$-group with $C_i$ of order $q^{f_i}$,
$m\le4$ and $f=\max_i f_i$. Given the $C_i$ explicitly and $R\in H$ implicitly, one finds
$\kappa_i\in[0,q^{f_i})$ with $R=\sum_i\kappa_iC_i$ using $O(f^2\log q)$ operations on implicit elements,
$O(f^2\log q+q^m)$ operations on explicit elements, and at most $fq^m$ comparisons between an implicit
and an explicit element.
\end{lemma}

\begin{proof}
Compute $R_j=q^jR$ for $0\le j<f$. We find, for $j=f,f-1,\dots,0$, integers
$\kappa_i^{(j)}\in[0,q^{\max(f_i-j,0)})$ with $q^jR=\sum_i\kappa_i^{(j)}q^jC_i$; for $j=f$ all are zero.
Given the $\kappa_i^{(j+1)}$, put $\Delta=R_j-\sum_i\kappa_i^{(j+1)}q^jC_i$. Then $q\Delta=0$ and
$\Delta\in q^jH=\bigoplus_{f_i>j}\langle q^jC_i\rangle$. The $q$-torsion of the cyclic group
$\langle q^jC_i\rangle$ is $\langle q^{f_i-1}C_i\rangle$, so $\Delta=\sum_{f_i>j}c_iq^{f_i-1}C_i$ for a unique
$(c_i)\in[0,q)^{m}$, which we find by comparing $\Delta$ with the at most $q^m$ explicit candidates. Put
$\kappa_i^{(j)}=\kappa_i^{(j+1)}+c_iq^{f_i-1-j}$ for $f_i>j$. Then $\kappa_i^{(j)}<q^{f_i-j}$ and
$\sum_i\kappa_i^{(j)}q^jC_i=(R_j-\Delta)+\Delta=R_j$. The candidate sets depend only on
$\{i:f_i>j\}$, so at most $m+1$ of them are built.
\end{proof}

Comparisons between an implicit and an explicit class are sound because reduced representations are unique.

For each prime $q\mid N$ with $q^e\,\|\,N$, put $N_q=N/q^e$, compute $D_q=N_qD$ implicitly and $C_{i,q}=N_qB_i$
explicitly. Since $N_q$ is prime to $q$, $C_{i,q}$ has order $q^{v_q(n_i)}$ and the $q$-Sylow subgroup of $J$
is $\bigoplus_i\langle C_{i,q}\rangle$. Lemma~\ref{lem:ph} gives $\kappa_{i,q}$ with
$D_q=\sum_i\kappa_{i,q}C_{i,q}$.

\subsection{Step 3: reading off the answer}
By the Chinese remainder theorem take $k_i\in[0,n_i)$ with $k_i\equiv\kappa_{i,q}\pmod{q^{v_q(n_i)}}$ for every
$q$, and compute $E=\sum_ik_iB_i$ explicitly. For every $q$ we have $N_q(D-E)=D_q-\sum_ik_iC_{i,q}=0$, so the order
of $D-E$ divides $\gcd_qN_q=1$ and $E=D$. By uniqueness, $E=(X-a',b)$; the algorithm reads off $a'$ and outputs
$a=a'-r$.

\subsection{Running time}
Step~1 costs $O(\ell^2)$ field operations in expectation. In Step~2, $\sum_q e_q\log q=\log N\le2\ell+1$, so all
implicit Jacobian arithmetic takes $O(\ell^2)$ group operations in $J$, hence $O(\ell^2)$ field operations; the
comparisons number at most $\sum_qe_qq^4\le(2\ell+1)y^4$, each needing $O(1)$ explicit constants. Step~3 is
explicit. Through Lemma~\ref{lem:dh-bbf} the total is $O(\ell^2y^4)$ group operations and $O(\ell^3)$ oracle
calls, that is $O(\ell^{4K+3})$. This proves Theorem~\ref{thm:main}.

\section{Fewer advice bits, and none}\label{sec:fewer}

The advice in Theorem~\ref{thm:main} is only a pointer to one good curve. Good curves are not rare, which
gives the three variants below.

\begin{lemma}\label{lem:density}
Let $K\ge K_0$ and $y=(\log x)^K$ as in Proposition~\ref{prop:goodcurve}. For large $p$, a uniformly random
coefficient vector $(a_0,\dots,a_5)\in\F_p^6$ defines a quintic $f=\sum a_iX^i$ with non-zero discriminant and
$y$-smooth $\#J_f(\F_p)$ with probability at least $p^{-2/K-o(1)}$.
\end{lemma}

\begin{proof}
At least $S_p-p\ge S_p/2$ of the $y$-smooth integers in the window avoid the exceptional set of
Theorem~\ref{thm:lpp}, and each is the order of at least $p^{9/2}/(24000\log(p^2)\log\log(p^2))$ quintics;
distinct orders give disjoint sets of quintics. Divide by $p^6$ and use $S_p\ge p^{3/2-2/K-o(1)}$.
\end{proof}

Without advice the algorithm computes $N=\#J_f(\F_p)$ by Pila's algorithm~\cite{Pil90,AH01} in time
$\ell^{O(1)}$, and a basis of $J_f(\F_p)$ in time $\ell^{O(1)}y^2$ from random elements
(Appendix~\ref{app:basis}).

\begin{theorem}[Polynomial time with $\varepsilon\ell$ yes/no queries]\label{thm:eps}
For every $0<\varepsilon\le1$ there is $c_\varepsilon$ such that, with at most $\varepsilon\ell$ queries to an oracle answering arbitrary yes/no questions, discrete
logarithms in groups of prime order $p$ are computed from a Diffie--Hellman oracle in expected time
$O(\ell^{c_\varepsilon})$ with failure probability at most $p^{-\varepsilon/2}$.
\end{theorem}

\begin{proof}
Take $K\ge\max(K_0,\lceil8/\varepsilon\rceil)$, so that a random coefficient vector is good with probability
$\lambda_p\ge p^{-\varepsilon/4-o(1)}$ by Lemma~\ref{lem:density}. Following~\cite[\S7]{KMS10} and Chor and
Goldreich~\cite{CG89}, identify $\F_p^6$ with $\F_{p^6}$, order it as $w_1,w_2,\dots$, draw $s,t\in\F_{p^6}$
and consider the pairwise independent samples $sw_j+t$, $j\le M=2^{\lfloor\varepsilon\ell\rfloor}$. By
Chebyshev's inequality none is good with probability at most $1/(M\lambda_p)\le p^{-\varepsilon/2}$ for large
$p$. The oracle supplies the $\lfloor\varepsilon\ell\rfloor$ bits of the least good index. The algorithm
then computes $N$ and a basis as above, and runs Steps 1--3 of Section~\ref{sec:reduction} in time
$\ell^{O(K)}$.
\end{proof}

Koblitz, Menezes and Shparlinski proved Theorem~\ref{thm:eps} under a conjecture on smooth numbers in short
intervals~\cite[Thm.~5]{KMS10}, and with time $L_p(\tfrac12)$ instead of polynomial under the Riemann
Hypothesis~\cite[Thm.~4]{KMS10}.

\begin{corollary}[One precomputation per prime]\label{cor:precomp}
For every $\varepsilon>0$ there is a uniform probabilistic algorithm in two phases. The first phase gets only
$p$, uses no group and no oracle, runs in expected time $p^{\varepsilon}\ell^{O(1)}$ and outputs a string
$\sigma_p$ of $O(\ell)$ bits, which anyone can check in time $\ell^{O(1/\varepsilon)}$. The second phase, given $\sigma_p$, any group $G$ of order $p$, a
Diffie--Hellman oracle for it and $g^a$, returns $a$ in expected time $\ell^{O(1/\varepsilon)}$.
\end{corollary}

\begin{proof}
Put $K=\max(K_0,\lceil4/\varepsilon\rceil)$. By Lemma~\ref{lem:density} a random coefficient vector is good with
probability at least $p^{-\varepsilon/2-o(1)}$. The first phase draws vectors, computes $\#J_f(\F_p)$ by Pila's
algorithm and tests $y$-smoothness by trial division, both in time $\ell^{O(1)}$ since $y\le(2\ell)^K$, until it
succeeds; after $p^\varepsilon$ failed draws it falls back to exhaustive search, which also covers the
finitely many small primes. It then computes a basis and a non-residue as above and outputs $\sigma_p$. The
second phase is Steps 1--3 of Section~\ref{sec:reduction}, in time $\ell^{4K+3}$.
\end{proof}

\begin{theorem}[No advice]\label{thm:uniform}
A probabilistic algorithm with a Diffie--Hellman oracle computes discrete logarithms in groups of prime order $p$ with unique element
encodings in expected time $L_p(\tfrac12,2)$, and in time $L_p(\tfrac12,\sqrt2)$ with batch smoothness
testing.
\end{theorem}

\begin{proof}
Take $y=\exp(c_0\sqrt{\log p\log\log p})$. Then $u\log u=(1/c_0+o(1))\sqrt{\log p\log\log p}$, and
Theorem~\ref{thm:smooth} together with Hildebrand's $\Psi(x,y)=x\rho(u)(1+o(1))$~\cite{Hil86} gives
$S_p\ge p^{3/2}\rho(u)^{1+o(1)}$, so a random coefficient vector is good with probability
$\rho(u)^{1+o(1)}=L_p(\tfrac12,1/c_0)^{-1}$. Draw vectors until one is good, testing smoothness of
$\#J_f(\F_p)$ by trial division in time $y\,\ell^{O(1)}$ per draw. Then compute $N$ and a basis in time
$\ell^{O(1)}y^2$ as above. Steps 1--3 with baby-step giant-step on each layer of Lemma~\ref{lem:ph} cost $q^2$
operations per layer, $(2\ell+1)y^2$ in all. With $c_0=1$ the total is $L_p(\tfrac12,2)$. A batch smoothness test~\cite{Ber04} makes each draw cost $\ell^{O(1)}$ amortised; then
the total is $L_p(\tfrac12,\max(1/c_0,2c_0))$, which is $L_p(\tfrac12,\sqrt2)$ at $c_0=1/\sqrt2$.
\end{proof}

\section{Consequences}\label{sec:consequences}

All reductions in this section are generic: in $G$ they use only the group operation, equality tests and the
oracle, and their other computations do not depend on how $G$ is encoded.

\subsection{Imperfect adversaries}
A security proof turns an adversary that solves Diffie--Hellman on some fraction of inputs into a
discrete-logarithm algorithm. Following Maurer and Wolf~\cite{MW96}, call $\mathcal O$ an
\emph{$\varepsilon$-Diffie--Hellman oracle} for $(G,g)$ if its output distribution on a given input does not
depend on earlier calls (it is time-invariant) and, on uniformly random input $(g^u,g^v)$, it returns $g^{uv}$
with probability at least $\varepsilon$, over the input and its own coins. A fixed
probabilistic algorithm run with fresh coins on each call is such an oracle.

\begin{corollary}\label{cor:imperfect}
Let $G=\langle g\rangle$ have prime order $p$ and let $\mathcal O$ be an $\varepsilon$-Diffie--Hellman oracle
for $(G,g)$. Given $\sigma_p$ and $\varepsilon$, discrete logarithms in $G$ are computed with an expected
$O(\ell^3\varepsilon^{-4}\log(\ell/\varepsilon))$ calls to $\mathcal O$ and $\mathrm{poly}(\ell,1/\varepsilon)$
further operations.
\end{corollary}

\begin{proof}
For every cyclic group of known order and every $\beta>0$, Maurer and Wolf~\cite[Thm.~1]{MW96} build from
$\mathcal O$ a Diffie--Hellman oracle that makes $O(\log(1/(\beta\varepsilon))/\varepsilon^4)$ calls to
$\mathcal O$ per query and answers each query correctly with probability at least $1-\beta$. The reduction of
Theorem~\ref{thm:main} makes $O(\ell^3)$ oracle calls in expectation; by Markov's inequality it makes at most
$T=O(\ell^3)$ calls with probability at least $5/6$, and each round is stopped after $T$ calls. With
$\beta=1/(6T)$ all answers in a round are correct with probability at least $5/6$, so a round succeeds with
probability at least $2/3$; the output $a$ is checked by computing $g^a$, and the round is repeated on
failure, at most $3/2$ times in expectation.
\end{proof}

\subsection{A characterisation for cyclic groups}
Let $n\ge2$ and let $q(n)$ be the largest prime $q$ with $q^2\mid n$, or $q(n)=1$ if $n$ is squarefree. Write
$\mathrm{GC}(n)$ for the least expected running time of a non-uniform generic reduction computing discrete
logarithms from a Diffie--Hellman oracle in cyclic groups of order $n$, in the sense of
Wolf~\cite[\S10.1]{Wol99}.

\begin{theorem}\label{thm:char}
$\sqrt{q(n)}/2-1\le\mathrm{GC}(n)\le(\log n)^{O(1)}\bigl(1+\sqrt{q(n)}\bigr)$. In particular discrete logarithm
and Diffie--Hellman are non-uniformly generically polynomial-time equivalent in cyclic groups of order $n$ if
and only if $q(n)\le(\log n)^{O(1)}$.
\end{theorem}

\begin{proof}
\emph{Upper bound.} The advice is the factorisation of $n$ and the strings $\sigma_q$ for the primes
$q\,\|\,n$, $O(\log n)$ bits in all. Pohlig--Hellman reduces discrete logarithms in $G$ to the subgroups of
prime-power order. For $q\,\|\,n$ put $h=g^{n/q}$; then $\DH_g(h^a,h^b)=h^{(n/q)ab}$, and raising to
$(n/q)^{-1}\bmod q$ gives a Diffie--Hellman oracle for $\langle h\rangle$, to which Theorem~\ref{thm:main}
applies with $\sigma_q$. For $q^e\,\|\,n$ with $e\ge2$, each of the $e$ base-$q$ digits costs $O(\sqrt q)$ by
baby-step giant-step.

\emph{Lower bound.} Maurer and Wolf~\cite{MW98}, \cite[Cor.~10.9]{Wol99}: every generic reduction of discrete
logarithms to a Diffie--Hellman oracle in cyclic groups of order $n$ has expected running time at least
$\sqrt q/2-1$, where $q$ is the largest prime with $q^2\mid n$.
\end{proof}

Wolf proved the equivalence statement~\cite[Thm.~10.10]{Wol99} under his Smoothness Assumption~1, which asks
for elliptic curves of polylogarithmically smooth order over every prime field.

\subsection{Black-box fields}
The algorithms of Sections~\ref{sec:reduction} and~\ref{sec:fewer} are $\F_p$-programs run on encodings
(Section~\ref{sec:prelim}), so they solve extraction for black-box fields isomorphic to $\F_p$ with the same
costs counted in field operations; Lemma~\ref{lem:dh-bbf} is only the translation back to groups. Maurer and
Raub~\cite[Lemma~2]{MR07} state this under ``a (plausible) number-theoretic conjecture'' and reduce the
problems for extension fields to it. They prove the consequences below for small characteristic, and for
large characteristic against non-uniform adversaries under the smoothness assumption
of~\cite{Mau94,MW99}.

\begin{corollary}\label{cor:bbf}
For every prime $p$ and every black-box field of characteristic $p$ and degree $k$ over $\F_p$, the
representation, extraction and isomorphism problems are solvable in time $\mathrm{poly}(k,\log p)$ with
$O(\log p)$ bits of advice depending only on $p$. Consequently, for finite fields whose characteristic is fixed
by the security parameter or chosen by the adversary, there is no field-homomorphic one-way permutation, and
no field-homomorphic encryption scheme that permits equality tests is one-way, against non-uniform
polynomial-time adversaries; in the private-key case plaintexts are recovered up to field isomorphism.
\end{corollary}

\begin{proof}
Maurer and Raub~\cite[Thm.~1, Cor.~2--3]{MR07} reduce these problems efficiently to the representation problem
for $\F_p$, which for $\F_p$ generated by $1$ is extraction; the black-box field form of
Theorem~\ref{thm:main} above solves extraction with $\sigma_p$. For the second part, such a scheme
instantiates a black-box field and extraction inverts it~\cite[Cor.~1, \S4]{MR07}.
\end{proof}

\section{Limitations and conclusion}\label{sec:conclusion}
The reduction is non-uniform: the advice names a curve for each characteristic, so the characteristic must
be fixed in advance, and finding a curve without advice costs $p^{2/K}$ trials, or subexponential time after
optimising $K$ (Theorem~\ref{thm:uniform}). A uniform polynomial-time reduction would need an efficient
construction of algebraic groups over $\F_p$ with polylogarithmically smooth order; Wolf conjectured that,
for general group orders, this is as hard as factoring~\cite[Conj.~10.1]{Wol99}. The constants are very
large: for $\delta=1/8$ the smoothness exponent is $K=25$ and the polynomial has degree about $4K$, so the
result says nothing about concrete security at cryptographic sizes.

Discrete logarithm and Diffie--Hellman are equivalent in every group of prime order, given $O(\log p)$ bits
that depend only on $p$. The obstacle was never the reduction but the distribution of smooth numbers, and a
zero-free half-plane is enough to supply them.

\section*{Use of AI}
The proof presented in this paper was developed through an interactive dialogue with Claude Fable~5.2 and
Claude Opus~5.5 (Anthropic). The author used them to explore multiple possible proof strategies and
directions, and the AI generated the initial draft under the author's guidance. The author then critically
examined the generated proof, verified each step for correctness and rigor, and revised it where necessary.
The author takes full responsibility for the final content.

\appendix
\section{Proof of Theorem~\ref{thm:smooth}}\label{app:smooth}

\paragraph{Notation.}
$L=\log x$, $S=\log y$, $u=L/S$, $Z_y(s)=\zeta(s,y)=\prod_{q\le y}(1-q^{-s})^{-1}$, $\alpha=\alpha(x,y)$ the
saddle point, $a_0=1-\alpha$, $W=x^\alpha Z_y(\alpha)$, $r=h/x$, and
$N(\sigma,T)=\#\{\rho=\beta+i\gamma:\zeta(\rho)=0,\ \beta\ge\sigma,\ |\gamma|\le T\}$.

\paragraph{Inputs from the literature.}
\begin{itemize}
\item[(I1)] $W\ll\Psi(x,y)\sqrt{LS}$ for $\log x<y\le x$~\cite[(1.3)]{You24}.
\item[(I2)] $a_0S=\log(u\log(u+1))+O(1)$ for $\log x<y\le x$~\cite[Lemma~2]{HT86}, \cite[(1.5)]{You24}. Hence
  there are absolute constants $0<c_*\le C_*$ with $c_*u\log(u+1)\le y^{a_0}\le C_*u\log(u+1)$ and
  $a_0S\le C_*\log(u+1)$, and $a_0S\ge1$ for $u$ large.
\item[(I3)] For fixed $\varrho>0$, $y\ge y_0(\varrho)$ with $\log x<y\le x$, and $1/S<|t|\le\exp(S^{3/2-\varrho})$,
  \[
  |Z_y(\alpha+it)|\le C_\varrho Z_y(\alpha)\exp\bigl(-c_0u\,t^2/(a_0^2+t^2)\bigr)
  \]
  \cite[Lemma~8]{HT86}, \cite[(1.6)]{You24}; we use $\varrho=1/10$.
\item[(I4)] For $2\le y,h,H\le x$~\cite[(1.9)]{You24},
  \[
  \Psi(x+h,y)-\Psi(x,y)=\frac1{2\pi i}\int_{\alpha-iH}^{\alpha+iH}Z_y(s)\frac{(x+h)^s-x^s}{s}\,ds+O(xL/H).
  \]
\item[(I5)] With $I_0$ the segment from $\alpha-i/S$ to $\alpha+i/S$~\cite[Lemma~12]{HT86},
  \cite[Lemma~3.6]{You24},
  \[
  \frac1{2\pi i}\int_{I_0}Z_y(s)x^s\,ds=\alpha\Psi(x,y)\Bigl(1+O\Bigl(\frac1u+\frac Sy\Bigr)\Bigr),
  \]
  and the same bound holds for $\frac1{2\pi}\int_{I_0}|Z_y(s)x^s|\,|ds|$.
\item[(I6)] Let $z\ge2$, $0<e<1/10$, $0<\eta\le1/2$, and suppose $\zeta$ has no zeros in
  $\{\operatorname{Re}s>1-\eta,\ |\operatorname{Im}s-t|<10z^\eta\}$. Then for $s=a+it\neq1$ with
  $1-(1-e)\eta\le a<1+1/\log z$,
  \begin{multline*}
  \Bigl|\sum_{n\le z}\frac{\Lambda(n)}{n^s}\Bigr|\le\frac{z^{1-a}+1}{|1-s|}+O_e\bigl((\log(|t|+2))^b\bigr)
  +O\bigl((\log z)^2z^{-e\eta}\bigr)\\
  +O_e\Bigl(\bigl((\log(|t|+z^\eta+2))^b+\tfrac1\eta\bigr)\log(z^\eta/\eta)\,z^{-e\eta/2}\Bigr),
  \end{multline*}
  with $b=b(e)<1$ and constants uniform in the other parameters~\cite[Lemma~2.3]{You24}.
\item[(I7)] For every $\mu>0$, $N(\sigma,T)\le C_\mu T^{30(1-\sigma)/13+\mu}$ for $\tfrac12\le\sigma\le1$,
  $T\ge2$~\cite[(1.4)]{GM24}.
\end{itemize}

\paragraph{Parameters.}
Fix $\delta,\theta$ and $K>K^*$. For small $\varepsilon,\varepsilon_1,\mu,\varepsilon_0>0$ put
$\tau=1-\theta+1/K+\varepsilon_1$, $g=(1-\varepsilon)^2-\varepsilon_0-\tfrac{30}{13}\tau$,
$\Xi_Z=\delta g-(1-\varepsilon)(1/K+\varepsilon_1)-\mu$ and
$\Xi_R=\theta-\tfrac12-\tfrac2K-2\varepsilon_1-\varepsilon$. At zero small parameters,
$\Xi_Z=[\delta(30\theta-17)-(13+30\delta)/K]/13>0$, which is $K>K^*$, and $\Xi_R=\theta-\tfrac12-2/K>0$ because
$K^*\ge56/(30\theta-17)>2/(\theta-\tfrac12)$. Fix $0<\varepsilon_0\le\varepsilon_1\le\varepsilon<1/12$ and
$\mu\le1/10$ with $\Xi_Z,\Xi_R>0$, and put $\kappa=\tfrac12\min(\Xi_Z,\Xi_R,\varepsilon_1/2)$. Then
$1/K+\varepsilon_1<(1-\varepsilon)\delta$. Let $x$ be large, $x^\theta\le h\le x$,
$(\log x)^K\le y\le x^{\varepsilon_0}$, and set $H=x^{1+1/K+\varepsilon_1}/h$, $T=H+10\sqrt y$,
$\sigma_{\mathrm R}=(1+\varepsilon)/2$, $D=D_0/L$ with $D_0=\max(8/\varepsilon,1/(c_*\log2))$, $t_0=Dy^{a_0}$ and
$t^*=Dy^{(1-\varepsilon)/2}$. For large $x$:
\begin{itemize}
\item[(P1)] $0<a_0\le1/K+\varepsilon_1/2$ ($\alpha$ increases with $y$, and (I2) at $y=(\log x)^K$);
\item[(P2)] $S\le\varepsilon_0L$, $x^{1/K+\varepsilon_1}\le H\le x^\tau$, $2\le T\le11x^\tau$;
\item[(P3)] $\Psi(x,y)\ge x^{1-a_0}/(CL)$ and $W\le C\Psi(x,y)L$;
\item[(P4)] $1\le D_0c_*\log(u+1)\le t_0S\le D_0C_*\log(u+1)$;
\item[(P5)] $a_0<(1-\varepsilon)\delta$, hence $\alpha>\sigma_{\mathrm R}$.
\end{itemize}

\begin{lemma}[local shift for every $y$]\label{lem:shift}
Fix $0<d\le\tfrac12$ and $0<e<1/10$. For $x\ge X_0(d,e)$, every $2\le y\le x$ and $d\le\eta\le\tfrac12$: let
$1-(1-e)\eta\le\sigma_1<\sigma_2<1+1/\log y$ and $0<|t|\le x$, suppose $\zeta$ has no zeros in
$\{\operatorname{Re}s>1-\eta,|\operatorname{Im}s-t|<10y^\eta\}$ and $|t|\ge\frac4{eL}(y^{1-\sigma_1}+1)$. Then
$|x^{\sigma_1+it}Z_y(\sigma_1+it)|\le x^{-(1-e)(\sigma_2-\sigma_1)}|x^{\sigma_2+it}Z_y(\sigma_2+it)|$.
\end{lemma}

\begin{proof}
For $\sigma_1\le a\le\sigma_2$ apply (I6) with $z=y$. The first term is at most
$(y^{1-\sigma_1}+1)/|t|\le eL/4$. The second is $O_e((2L)^b)$. The third is
$(\log y)^2y^{-e\eta}=w^2e^{-e\eta w}$ with $w=\log y$, at most $4/(e\eta\,\mathrm e)^2\le4/(ed)^2$ for every
$y\ge2$. In the fourth, $\log(|t|+y^\eta+2)\le2L$, $1/\eta\le1/d$ and
$\log(y^\eta/\eta)y^{-e\eta/2}\le2/(e\,\mathrm e)+\log(1/d)$, so it is $O_{d,e}(L^b+1)$. As $b<1$, the
total is at most $eL/2$ for large $x$. Since $-Z_y'/Z_y(s)=\sum_{n\le y}\Lambda(n)n^{-s}+R_y(s)$ with
$|R_y(s)|\le\sum_q(\log q)q^{-2a}/(1-q^{-a})=O_e(1)$ for $a\ge(1+e)/2$, we get $|Z_y'/Z_y(a+it)|\le eL$, and
integrating $\partial_a\log|Z_y(a+it)|$ from $\sigma_1$ to $\sigma_2$ gives the claim.
\end{proof}

\paragraph{The contour.}
Let $\mathcal Z_T$ be the zeros $\rho=1-\nu+i\gamma$ with $\nu\le\tfrac12$ and $|\gamma|\le T$; by hypothesis
$\nu\ge\delta$. For $t\in\mathbb R$ let $\sigma_0(t)=\alpha$ for $|t|\le t_0$, $\sigma_0(t)=1-\log(|t|/D)/S$ for
$t_0\le|t|\le t^*$ and $\sigma_0(t)=\sigma_{\mathrm R}$ beyond; and for $\rho\in\mathcal Z_T$ let
$\sigma_\rho(t)=1-(1-\varepsilon)\nu$ for $|t-\gamma|\le10y^\nu$,
$\sigma_\rho(t)=1-(1-\varepsilon)\log(|t-\gamma|/10)/S$ for $10y^\nu\le|t-\gamma|\le10\sqrt y$ and
$\sigma_\rho(t)=\sigma_{\mathrm R}$ beyond. Put $\sigma_\Gamma=\max(\sigma_0,\max_\rho\sigma_\rho)$,
$\Gamma=\{\sigma_\Gamma(t)+it:|t|\le H\}$, and let $C_\pm$ be the horizontal segments from $\Gamma$ to
$\alpha\pm iH$. By (P5), $\sigma_\rho\le1-(1-\varepsilon)\delta<\alpha$, so $\sigma_\Gamma\le\alpha$ with
equality exactly for $|t|\le t_0$; and $t_0\ge1/S$ by (P4).

\begin{lemma}\label{lem:pointwise}
For $|t|\le H$ and $\sigma_\Gamma(t)\le\sigma<\alpha$:
\[
|x^{\sigma+it}Z_y(\sigma+it)|\le x^{-(1-\varepsilon)(\alpha-\sigma)}\,|x^{\alpha+it}Z_y(\alpha+it)|
\le x^{-(1-\varepsilon)(\alpha-\sigma)}W .
\]
\end{lemma}

\begin{proof}
Let $\eta=\max(\delta,(1-\sigma)/(1-\varepsilon))\le\tfrac12$. No zero lies in
$\{\operatorname{Re}s>1-\eta,|\operatorname{Im}s-t|<10y^\eta\}$: for $\eta=\delta$ this is the hypothesis;
otherwise such a zero would be in $\mathcal Z_T$ and would force $\sigma_\rho(t)>\sigma$. Next,
$\sigma<\alpha$ forces $|t|>t_0$, and then $Dy^{1-\sigma}\le|t|$ from the definition of $\sigma_0$ (using
$\sigma\ge\sigma_{\mathrm R}$ beyond $t^*$), so $|t|\ge\frac8{\varepsilon L}y^{1-\sigma}\ge\frac4{\varepsilon L}(y^{1-\sigma}+1)$.
Lemma~\ref{lem:shift} with $d=\delta$, $e=\varepsilon$, $\sigma_1=\sigma$, $\sigma_2=\alpha$ gives the first
inequality, and $|Z_y(\alpha+it)|\le Z_y(\alpha)$ the second.
\end{proof}

For $\operatorname{Re}s\le1$ one has $|((x+h)^s-x^s)/s|\le r|x^s|$ and
$|((x+h)^s-x^s)/s-rx^s|\le\tfrac12r^2|s-1||x^s|$; on $C_\pm$, $|((x+h)^s-x^s)/s|\le3|x^s|/H$.

\begin{lemma}[zeros]\label{lem:zeros}
$\sum_{\rho\in\mathcal Z_T}\int_{\Gamma\cap\Gamma_\rho}|x^sZ_y(s)||ds|\le CWx^{-\Xi_Z}\le C\Psi(x,y)x^{-\kappa}$,
where $\Gamma_\rho$ is the part of the contour defined by $\sigma_\rho$.
\end{lemma}

\begin{proof}
Let $c_1=(1-\varepsilon)^2-S/L\ge(1-\varepsilon)^2-\varepsilon_0$. On the vertical piece (length $20y^\nu$),
Lemma~\ref{lem:pointwise} gives at most $20Wx^{(1-\varepsilon)a_0}x^{-c_1\nu}$. On the curved pieces,
parametrised by $v\in[\nu,\tfrac12]$ with $|s'(v)|\le11Sy^v$, the integral is at most
$11(S/((1-\varepsilon)^2L-S))Wx^{(1-\varepsilon)a_0}x^{-c_1\nu}\le11Wx^{(1-\varepsilon)a_0}x^{-c_1\nu}$. Group the
zeros by $m=\lfloor\nu L\rfloor\ge\lfloor\delta L\rfloor=m_0$. By (I7) and (P2) there are at most
$121C_\mu x^\mu e^{(30/13)\tau(m+1)}$ zeros with this $m$, each contributing at most $42Wx^{(1-\varepsilon)a_0}e^{-c_1m}$.
Since $c_1-\tfrac{30}{13}\tau\ge g>0$, the sum over $m\ge m_0$ is
$\le C(\mu,g)Wx^{(1-\varepsilon)a_0+\mu-g\delta}\le CWx^{-\Xi_Z}$ by (P1); finally $W\le C\Psi L$.
\end{proof}

\begin{lemma}[remainder line and connectors]\label{lem:remainder}
The part of $\Gamma$ on $\operatorname{Re}s=\sigma_{\mathrm R}$ contributes at most $2Wx^{-\Xi_R}\le C\Psi x^{-\kappa}$, and
$\int_{C_\pm}|x^sZ_y(s)||ds|\le W/((1-\varepsilon)L)\le C\Psi/\sqrt u$.
\end{lemma}

\begin{proof}
On $\operatorname{Re}s=\sigma_{\mathrm R}$, Lemma~\ref{lem:pointwise} gives
$|x^sZ_y|\le Wx^{-(1-\varepsilon)(\alpha-\sigma_{\mathrm R})}$ with $(1-\varepsilon)(\alpha-\sigma_{\mathrm R})\ge\tfrac12-\varepsilon-1/K-\varepsilon_1$, over length $\le2x^\tau$; the exponent is
$-\Xi_R$. On $C_\pm$, $\int_{-\infty}^\alpha x^{-(1-\varepsilon)(\alpha-\sigma)}d\sigma=1/((1-\varepsilon)L)$, and
$W/L\le C\Psi\sqrt{S/L}$ by (I1).
\end{proof}

\begin{lemma}[central tails and curves]\label{lem:central}
Let $\Gamma_0^\pm$ be the exponential paths $\{1-v\pm iDy^v:a_0\le v\le(1-\varepsilon)/2\}$. Then
\[
\int_{\substack{s=\alpha+it\\1/S<|t|\le t_0}}|x^sZ_y(s)|\,|ds|+\int_{\Gamma\cap\Gamma_0^\pm}|x^sZ_y(s)|\,|ds|
\le C\,\frac{\Psi(x,y)}u .
\]
\end{lemma}

\begin{proof}
For $1/S<|t|\le t_0$, (I3) applies and $t^2/(a_0^2+t^2)\ge c/(\log(u+1))^2$ by (I2), so with (P4) and (I1)
the tail is at most $CD_0\Psi\sqrt u\log(u+1)e^{-cu/\log^2(u+1)}\le C'\Psi/u$. On $\Gamma_0^\pm$, $|t|\ge t_0$
gives $t^2/(a_0^2+t^2)\ge c$ and (I3) gives $|Z_y(\alpha+it)|\le Ce^{-cu}Z_y(\alpha)$; combining this with the first
inequality of Lemma~\ref{lem:pointwise} and $|s'(v)|\le2DSy^v$ the curves contribute at most
$CWe^{-cu}D_0\log(u+1)/L\le C'\Psi/u$.
\end{proof}

\paragraph{Assembly.}
Let $\Gamma'=\Gamma\setminus I_0$. By (I4) and Cauchy's theorem,
\[
\Psi(x+h,y)-\Psi(x,y)=\frac1{2\pi i}\Bigl(\int_{I_0}+\int_{\Gamma'}+\int_{C_-}+\int_{C_+}\Bigr)
Z_y(s)\frac{(x+h)^s-x^s}{s}\,ds+O\Bigl(\frac{xL}H\Bigr).
\]
By (I5) in both forms and $\max|s-1|\le a_0+1/S$ on the central
segment, the central part is $\alpha r\Psi(1+O(1/u))+\vartheta\tfrac12r^2(a_0+1/S)\alpha\Psi(1+O(1/u))$ with
$|\vartheta|\le1$. Lemmas~\ref{lem:zeros}--\ref{lem:central} bound the rest of $\Gamma$ by
$Cr\Psi(1/u+x^{-\kappa})$; the connectors give $Cr\Psi x^{-1/K-\varepsilon_1}/\sqrt u$ since $1/H=rx^{-1/K-\varepsilon_1}$;
and $xL/H=rx^{1-1/K-\varepsilon_1}L\le Cr\Psi L^2x^{a_0-1/K-\varepsilon_1}\le Cr\Psi x^{-\kappa}$ by (P1), (P3).
Dividing by $\alpha r\Psi$ (with $\alpha\ge\tfrac12$ and $x^{-\kappa}\le1/u$),
\[
\Bigl|\frac{\Psi(x+h,y)-\Psi(x,y)}{\alpha r\Psi(x,y)}-1\Bigr|\le\frac{C_1}u+\frac r2(a_0+1/S)\Bigl(1+\frac{C_1}u\Bigr),
\]
which gives \eqref{eq:asym} since $a_0+1/S\le(C_*+2)\log(u+1)/S$. For the lower bound,
$a_0+1/S\le1/K+\varepsilon_1<\tfrac12$, so the right side is at most $\tfrac12$ once $u\ge5C_1$, which holds after
decreasing $\varepsilon_0$; and $\alpha\ge1-1/K-\varepsilon_1\ge\tfrac23$. \qed

\paragraph{Constants.}
For $\theta=3/4$ and $\delta=1/8$ one gets $K^*=268/11<25$, the value $K_0=25$ of
Proposition~\ref{prop:goodcurve}. The half-planes $\operatorname{Re}s>11/12$~\cite{Kon26} and
$\operatorname{Re}s>31/32$~\cite{CT26} give $K_0=34$ and $K_0=82$. With a zero-density exponent $A$ in place of
Guth and Maynard's $30/13$, the threshold is $K^*_A=(1+A\delta)/(\delta(1-A(1-\theta)))$.

\paragraph{Under the Riemann Hypothesis.}
All zeros then lie on $\operatorname{Re}s=\tfrac12$, every $\Gamma_\rho$ lies on the remainder line,
Lemma~\ref{lem:zeros} is void, and only $\Xi_R>0$, (P5) and the Perron step are used: $K>2/(\theta-\tfrac12)$
suffices, that is $K=9$ for $\theta=3/4$, which makes the constant of Younis~\cite[Thm.~1.3]{You24} explicit in
this range.

\section{Computing the group without advice}\label{app:basis}
Let $f$ be a monic quintic with non-zero discriminant such that $N=\#J_f(\F_p)$ is $y$-smooth, and factor $N$
by trial division up to $y$. Draw a monic quadratic $u$ uniformly, let $m(u)\in\{0,1,2,4\}$ be the number of $v$ with
$\deg v\le1$ and $v^2\equiv f\pmod u$ (computed with square roots in $\F_p$ or $\F_{p^2}$), accept with
probability $m(u)/4$ and output a uniform one of these $v$. Every class with $\deg u=2$ is output with the
same probability, and the remaining classes, at most $p+4\sqrt p+1$ of them, form an $O(1/p)$ fraction of
$J_f(\F_p)$. Six such samples generate $J_f(\F_p)$ with probability at least $0.29-O(1/p)$, by Pomerance's
bound for groups of rank at most four~\cite{Pom02} as used in~\cite[\S7]{KMS10}. Teske's
algorithm~\cite[Thm.~6.1]{Tes99}, in the form of Sutherland~\cite[\S1]{Sut11}, turns them into a basis: for a
$q$-group of rank $m\le4$ it costs $\ell^{O(1)}q^{m/2}\le\ell^{O(1)}q^2$ group operations, so $\ell^{O(1)}y^2$ in
all; the algorithm checks that the basis has order $N$, and draws again otherwise.

\bibliographystyle{alpha}
\bibliography{refs}
\end{document}